\documentclass[11pt,reqno]{amsart}
\usepackage{amsmath,amssymb,geometry,color,tikz-cd}
\usepackage{comment}
\usepackage[backref,pagebackref,pdftex,hyperindex]{hyperref}
\numberwithin{equation}{section}%
\usepackage[alphabetic]{amsrefs}
\usepackage{amssymb}
\usepackage{bbm}
\usepackage{enumitem}
\usepackage{upgreek}

\newtheorem{proposition}{Proposition}[section]
\newtheorem{theorem}[proposition]{Theorem}
\newtheorem{lemma}[proposition]{Lemma}

\theoremstyle{definition}
\newtheorem{remark}[proposition]{Remark}
\newtheorem{definition}[proposition]{Definition}

\newtheorem{question}[proposition]{Question}
\newtheorem{claim}[proposition]{Claim}

\title{Degeneration of log Calabi-Yau pairs via log canonical places}
\author{Guodu Chen and Chuyu Zhou}

\address{Institute for Theoretical Sciences, Westlake University, Hangzhou, Zhejiang, 310024, China}
\email{chenguodu@westlake.edu.cn}
	
\address{\'Ecole Polytechnique F\'ed\'erale de Lausanne (EPFL), MA C3 615, Station 8, 1015 Lausanne, Switzerland}
\email{chuyu.zhou@epfl.ch}
\date{} 

\thanks{2010 
	    \emph{Mathematics Subject Classification}: 14E30, 14J17.
	    \newline
	    \indent 
		\emph{Keywords}: log Calabi-Yau pair, log canonical place, weakly special test configuration.
	}

\newcommand{\ord}{{\rm {ord}}}
\newcommand{\tc}{{\rm {tc}}}

\newcommand{\red}{{\rm {red}}}
\newcommand{\Spec}{{\rm {Spec}}}
\newcommand{\Proj}{{\rm{Proj}}}

\newcommand{\Supp}{{\rm {Supp}}}

\usepackage{todonotes}

\newcommand{\bA}{\mathbb{A}}

\newcommand{\bC}{\mathbb{C}}

\newcommand{\bN}{\mathbb{N}}

\newcommand{\bQ}{\mathbb{Q}}
\newcommand{\bR}{\mathbb{R}}

\newcommand{\bZ}{\mathbb{Z}}

\newcommand{\mB}{\mathcal{B}}

\newcommand{\mD}{\mathcal{D}}
\newcommand{\mE}{\mathcal{E}}
\newcommand{\mF}{\mathcal{F}}

\newcommand{\mL}{\mathcal{L}}

\newcommand{\mO}{\mathcal{O}}

\newcommand{\mR}{\mathcal{R}}

\newcommand{\mW}{\mathcal{W}}
\newcommand{\mX}{\mathcal{X}}
\newcommand{\mY}{\mathcal{Y}}
\newcommand{\mZ}{\mathcal{Z}}

\newcommand{\tZ}{\tilde{Z}}

\newcommand{\ka}{\mathfrak{a}}

\begin{document}

\begin{abstract}
Let $(X,\Delta)$ be a projective log canonical Calabi-Yau pair and $L$ an ample $\bQ$-line bundle on $X$, we show that there is a correspondence between lc places of $(X,\Delta)$ and weakly special test configurations of $(X,\Delta; L)$. 
\end{abstract}

\maketitle
\tableofcontents

\section{Introduction}

We work over $\bC$ throughout. 

Test configurations are basic objects in the study of K-stability of Fano varieties. After the work \cite{BHJ17}, people realize that the language via test configurations can be reformulated by a even more basic concept, i.e. valuations. As a special example of valuations, divisorial valuations play important roles in the recent development of algebraic K-stability theory of Fano varieties (e.g. \cite{Xu21}). For example, any weakly special test configuration with integral central fiber of a $\bQ$-Fano variety (i.e. a Fano variety with klt singularities) could be induced by an lc place of some complement (see \cite[Appendix]{BLX19} or \cite{CZ21}), and via reformulating these test configurations by the corresponding lc places of complements, many tools in birational geometry come to play the powerful role. In the realm of Calabi-Yau varieties, K-stability theory could be translated into singularity theory by \cite{Oda13} (see also \cite[Section 9]{BHJ17}), in other words, people are mainly interested in Calabi-Yau varieties with log canonical or semi-log canonical singularities. Inspired by this, we aim to reformulate  test configurations via valuations (more precisely via divisorial valuations induced by lc places) for Calabi-Yau varieties, and expect to further study the degenerations of Calabi-Yau varieties in the moduli theory.

We say that a projective log pair $(X,\Delta)$ is a \textit{log canonical} (resp. \textit{semi-log canonical}) \textit{Calabi-Yau pair} (abbreviated by \textit{lc CY pair} (resp. \textit{slc CY pair})) if $(X,\Delta)$ admits log canonical (resp. semi-log canonical) singularities and $K_X+\Delta\sim_\bQ 0$. We aim to establish the following result.

\begin{theorem}\label{thm: main1}
Let $(X,\Delta)$ be a projective lc CY pair and $L$ an ample $\bQ$-line bundle on $X$, then we have the following result:
\begin{enumerate}
\item If $E$ is an lc place of $(X,\Delta)$, then $E$ is dreamy with respect to $L$;
\item there is a correspondence between lc places of $(X,\Delta)$ and non-trivial weakly special test configurations of $(X,\Delta;L)$  with integral central fibers. More precisely, given an lc place $E$ of $(X,\Delta)$ and a positive integer $c\in \bZ^+$, there exists a weakly special test configuration with integral central fiber $(\mX,\Delta_\tc; \mL)\to \bA^1$ of $(X,\Delta; L)$ such that ${\ord_{\mX_0}}|_{K(X)}=c\cdot \ord_E$; conversely, given a weakly special test configuration with integral central fiber $(\mX,\Delta_\tc; \mL)\to \bA^1$ of $(X,\Delta; L)$ such that ${\ord_{\mX_0}}|_{K(X)}=c\cdot \ord_E$ for some $c\in \bZ^+$, then $E$ is an lc place of $(X,\Delta)$.
\end{enumerate}
\end{theorem}

We just note here that $E$ is dreamy with respect to $L$ means the finite generation of the following graded ring
$$\bigoplus_{m\in \bN}\bigoplus_{j\in \bN} H^0(Y, g^*(mrL)-jE), $$
where $g: Y\to X$ is a normal projective birational model such that $E$ is a prime divisor on $Y$, and $r$ is a fixed positive integer such that $rL$ is Cartier. We say that the test configuration $(\mX,\Delta_\tc; \mL)\to \bA^1$ of $(X,\Delta; L)$ is weakly special if $(\mX, \Delta_\tc+\mX_0)$ is log canonical (see Definition \ref{defn: tc}). Theorem \ref{thm: main1} concerns weakly special test configuration with integral central fiber, while in Section \ref{sec: g-global correspondence}, we study the case where the central fiber is reduced. The key ingredient to establish the above global correspondence is the following characterization of log canonical blow-ups (see Definition \ref{def: lc blowup}).

\begin{theorem}{\rm{(= Theorem \ref{thm: key})}}
Let $(X, \Delta)$ be a projective lc CY pair and $L$ an ample $\bQ$-line bundle on $X$. Denote by $Z=\Spec\ \bigoplus_{m\in \bN} H^0(X, mrL)$, 
where $r$ is a positive integer such that $rL$ is Cartier. Let $\Delta_Z$ be the extension of $\Delta$ on $Z$ and $o\in Z$ the cone vertex of $Z$.
Suppose $w$ is a $\bC^*$-equivariant divisorial valuation corresponding to an lc blow-up over $o\in Z$, then $w$ is a  quasi-monomial combination of $v_0$ and $\ord_{E_\infty}$ with weight $(\lambda, c)$ for some $\lambda\in \bZ^{+}$ and $c\in \bQ^{\geq 0}$,
where $v_0$ is the canonical valuation, $E_\infty$ is the natural extension of some prime divisor $E$ over $X$, and $E$ is an lc place of $(X,\Delta)$.
\end{theorem}

In Section \ref{sec: g-global correspondence}, we will generalize Theorem \ref{thm: main1} for weakly special test configurations with reduced central fibers, and establish a correspondence between weakly special collections (see Definition \ref{g-lc blowup}) and weakly special semi-test configurations (see Definition \ref{defn: tc}).

\begin{theorem}{\rm{(=Theorem \ref{thm: g-global correspondence})}}\label{thm: main3}
Let $(X, \Delta)$ be a projective lc CY pair and $L$ an ample $\bQ$-line bundle on $X$, then there exists a  correspondence between weakly special collections of $(X, \Delta)$ and weakly special semi-test configurations of $(X, \Delta; L)$. More precisely, given a weakly special collection $\{c_1\cdot \ord_{E_1},...,c_l\cdot \ord_{E_l}\},$ we can construct a weakly special semi-test configuration  
$(\mX, \Delta_\tc; \mL)$ satisfying 
$${\ord_{\mX_{0,i}}}|_{K(X)}=c_i\cdot \ord_{E_i}, i=1,...,l,$$
where $\mX_{0,i}$ are all non-trivial components of $\mX_0$; conversely, given a weakly special semi-test configuration  
$(\mX, \Delta_\tc; \mL)$ satisfying 
$${\ord_{\mX_{0,i}}}|_{K(X)}=c_i\cdot \ord_{E_i}, i=1,...,l,$$
where $c_i\in \bZ^+$ and $\mX_{0,i}$ are all non-trivial components of $\mX_0$, then the set $\{c_1\cdot \ord_{E_1},...,c_l\cdot \ord_{E_l}\}$ is a weakly special collection.
\end{theorem}

As a corollary, we see that there are no non-trivial weakly special semi-test configurations of $(X,\Delta;L)$ if $(X,\Delta)$ is a klt Calabi-Yau pair (e.g. \cite[Corollary 4.3]{Odaka12}). 

We emphasize here that our definition of weakly special collection is a little different from that in \cite[Appendix A]{BLX19}, as they define weakly special collections directly from  weakly special semi-test configurations, while our definition automatically contains the information of lc places (see Definition \ref{g-lc blowup}). Thus the above theorem is also a correspondence between weakly special semi-test configurations and lc places of a given polarized lc CY pair.

\begin{remark}
For Theorem \ref{thm: main1} and Theorem \ref{thm: main3}, we have explained precisely what the correspondences mean in the both statements. We emphasize that such correspondences are not bijective. For example, given a weakly special collection as in Theorem \ref{thm: main3}, we could construct many weakly special semi-test configurations of $(X, \Delta; L)$ satisfying the requirements (see the proof of Theorem \ref{thm: g-global correspondence} and Remark \ref{rem: set of tc}).
\end{remark}

\
\

\textit{Notation.} For a given normal variety $X$, we say that $\Delta$ is a $\bQ$-divisor on $X$ if $\Delta$ can be put as a finite sum $\Delta:=\sum a_i \Delta_i$, where $a_i\in \bQ$ and $\Delta_i$ are Weil divisors on $X$. We say that $\Delta$ is $\bQ$-Cartier if $m\Delta$ is a Cartier divisor for a sufficiently divisible integer $m\in \bZ^+$. We say that $L$ is a $\bQ$-line bundle on $X$ if $L^{\otimes m}$ (or put as $mL$) is a line bundle on $X$ for a sufficiently divisible integer $m\in \bZ^+$.

We say that $(X, \Delta)$ is a log pair if $X$ is a normal variety and $\Delta$ is an effective $\bQ$-divisor on $X$ such that $K_X+\Delta$ is $\bQ$-Cartier.
We say that $E$ is a prime divisor over $X$ if there is a proper birational morphism from a normal variety $Y$, denoted by $Y\to X$, such that $E$ is a prime divisor on $Y$. 

For a given normal variety $X$, we say that $v$ is a divisorial valuation over $X$ if it is of the form $v=c\cdot \ord_E$, where $c\in \bR^+$ and $E$ is a prime divisor over $X$. We refer to \cite[Section 1.3]{BHJ17} for a basic introduction of valuation.

Given a normal projective variety $X$, suppose $L$ is an ample line bundle on $X$ and $Z:=\Spec\ \oplus_{m\in \bN} H^0(X, mL)$ is the affine cone over $X$ with respect to $L$. Let $o\in Z$ be the cone vertex, then $Z\setminus o$ is a $\bC^*$-bundle over $X$. Thus any prime divisor over $X$ admits a natural $\bC^*$-equivariant extension over $Z$. More generally, any valuation over $X$ can be naturally extend to a $\bC^*$-invariant valuation over $Z$.  Suppose $\tilde{Z}\to Z$ is the blow-up of $o\in Z$, we write $v_0$ for the divisorial valuation corresponding to the exceptional divisor which is isomorphic to $X$ ($v_0$ is also called canonical valuation). Let $v$ be a valuation over $X$ and $\bar{v}$ the valuation over $Z$ via the natural extension, then we put the quasi-monomial combination of $v_0$ and $\bar{v}$ with weight $(\lambda, c)$ as $\lambda\cdot v_0+c\cdot \bar{v}$ for convenience\footnote{We emphasize here that the valuation $\lambda\cdot v_0+c\cdot \bar{v}$ does not necessarily have the property $(\lambda\cdot v_0+c\cdot \bar{v})(f)=\lambda\cdot v_0(f)+c\cdot \bar{v}(f)$ for every $f\in K(Z)^*$, it is just a notation of quasi-monomial combination for convenience.}.

For various type of singularities in birational geometry such as lc, klt, dlt, slc, etc., we refer to \cite{KM98, Kollar13}.

\

\
\noindent
{\bf Acknowledgements}:
G. Chen is supported by the China post-doctoral grants BX2021269 and 2021M702925. C. Zhou is supported by grant European Research Council (ERC-804334). We would like to thank Yuchen Liu for helpful comments and suggestions, and we also thank the referees for careful reading and providing a list of precious comments.

\

\section{Local correspondence}\label{sec: local correspondence}

In this section, we study weakly special test configuration from a local point of view. We first recall some basic concepts.

\begin{definition}
Let $(X,\Delta)$ be a log pair and $E$ a prime divisor over $X$, then we define the \emph{log discrepancy} of $E$ with respect to $(X,\Delta)$ as follows:
$$A_{X,\Delta}(E)=\ord_E(K_Y-g^*(K_X+\Delta))+1, $$
where $g: Y\to X$ is a log resolution such that $E$ is a prime divisor on $Y$. If $A_{X, \Delta}(E)=0$, we say that $E$ is an lc place of $(X, \Delta)$.
\end{definition}

\begin{definition}\label{defn: tc}
Let $(X,\Delta)$ be a projective log pair and $L$ an ample $\bQ$-line bundle on $X$. A \emph{semi-test configuration} $\pi: (\mX,\Delta_\tc;\mL)\to \bA^1$ is a family over $\bA^1$ consisting of the following data:
\begin{enumerate}
\item $\pi: \mX\to \bA^1$ is a projective flat morphism from a normal variety $\mX$, $\Delta_\tc$ is an effective $\bQ$-divisor on $\mX$ and $\mL$ is a relatively semi-ample $\bQ$-line bundle on $\mX$,
\item the family $\pi$ admits a $\bC^*$-action which lifts the natural $\bC^*$-action on $\bA^1$ such that $(\mX,\Delta_\tc; \mL)\times_{\bA^1}\bC^*$ is $\bC^*$-equivariantly isomorphic to $(X, \Delta; L)\times_{\bA^1}\bC^*$.
\end{enumerate}
We say that the semi-test configuration is a \emph{test configuration} if $\mL$ is relatively ample. We say that the (semi-)test configuration is \emph{weakly special} if $(\mX, \Delta_\tc+\mX_0)$ is log canonical. In this case, we also say that $(\mX, \Delta_\tc)$ is a weakly special degeneration of $(X, \Delta)$. For a 
given weakly special (semi-)test configuration $(\mX, \Delta_\tc;\mL)$, if the central fiber $\mX_0$ is integral, we say that it is a \textit{weakly special (semi-)test configuration with integral central fiber}.

\end{definition}

\begin{definition}\label{def: lc blowup}
Let $(X,\Delta)$ be a log pair and $x\in X$ is a closed point. Let $E$ be a prime divisor over $X$, we say that $E$ is a \textit{log canonical blow-up} (abbreviated by \textit{lc blow-up}) over $x$ if there exists a proper birational morphism $f: Y\to X$ such that $Y\setminus E\cong X\setminus x$, $(Y, f_*^{-1}{\Delta}+E)$ is log canonical, and $-E$ is relatively ample.
\end{definition}

In the rest of this section, we fix the following notation.

Let $(X,\Delta)$ be a projective lc CY pair and $L$ an ample $\bQ$-line bundle on $X$. Consider the affine cone over $X$ with respect to $rL$:
$$Z:=\Spec\ \bigoplus_{m\in \bN} H^0(X, mrL), $$
where $r$ is a positive integer such that $rL$ is Cartier. Denote by $o\in Z$ the cone vertex and $\Delta_Z$ the extension of $\Delta$ on $Z$. We need the following theorem to establish the global correspondence in Section \ref{sec: global correspondence}.

\begin{theorem}\label{thm: key}
Suppose $w$ is a $\bC^*$-equivariant divisorial valuation corresponding to an lc blow-up over $o\in Z$, then $w$ is a  quasi-monomial combination of $v_0$ and $\ord_{E_\infty}$ with weight $(\lambda, c)$ for some $\lambda\in \bZ^+, c\in \bQ^{\geq 0}$,
where $E_\infty$ is the natural extension of some prime divisor $E$ over $X$, and $E$ is an lc place of $(X,\Delta)$.
\end{theorem}

\begin{proof}
If $w$ is obtained by the blow-up of $o\in Z$: $\tilde{Z}\to Z$, then $\lambda=1, c=0$. We may assume that $w$ is not obtained via this blow-up. By \cite[Lemma 4.2]{BHJ17}, we can express $w$ as a  quasi-monomial combination of $v_0$ and $\overline{r(w)}$ with weight $(\lambda, 1)$ for some $\lambda\in \bZ^+$,
$$w\ =\ \lambda\cdot v_0+\overline{r(w)}, $$
where $r(w)=w|_{K(X)}$, and $\overline{r(w)}$ is the natural extension of $r(w)$. As $w$ is a divisorial valuation, by \cite[Lemma 4.1]{BHJ17}, $r(w)$ is also divisorial. Thus one can find some $c\in \bR^+$ and a prime divisor $E$ over $X$ such that $r(w)=c\cdot \ord_E$. Denote by $E_\infty$ the natural extension of $E$, then $w$ is a  quasi-monomial combination of $v_0$ and $\ord_{E_\infty}$ with weight $(\lambda, c)$,
$$w\ =\ \lambda\cdot v_0+c\cdot \ord_{E_\infty}. $$
As $w$ is a divisorial valuation, the value group of $w$ has rational rank $1$, which implies that $c\in \bQ$.
It remains to show that $E$ is an lc place of $(X,\Delta)$.
We denote $f: W\to Z$ to be the lc blow-up which gives rise to $w$, and $G\subset W$ the exceptional divisor. Then we have $w=\ord_G$ and
$$K_W+f_*^{-1}\Delta_Z+G=f^*(K_Z+\Delta_Z)+A_{Z,\Delta_Z}(G) G,$$
where $(W, f_*^{-1}\Delta_Z+G)$ is log canonical. If $A_{Z,\Delta_Z}(G)=0$, then $G$ is already an lc place of $(Z,\Delta_Z)$. As $v_0$ is an lc place of $(Z,\Delta_Z)$ (see \cite[Proposition 3.14]{Kollar13}), by \cite[Proposition 5.1]{JM12}, $\overline{r(w)}$ is also an lc place of $(Z,\Delta_Z)$, which implies that $E$ is an lc place of $(X,\Delta)$. From now on, we assume $A_{Z,\Delta_Z}(G)>0$ and derive the contradiction. 

Denote by
$$K_G+\Delta_G:=(K_W+f_*^{-1}\Delta_Z+G)|_G, $$
then $(G, \Delta_G)$ is semi-log canonical by the sub-adjunction, and 
$$-(K_G+\Delta_G)\sim_\bQ -A_{Z,\Delta_Z}(G)G|_G$$ 
is ample. Choose a sufficiently divisible $k\in \bZ^+$ such that $$k(K_W+f_*^{-1}\Delta_Z+G),\  kG \text{ and } k(K_Z+\Delta_Z)$$ 
are all Cartier. Consider the following exact sequence:
$$0\to \mO_W(-G-k(K_W+f_*^{-1}\Delta_Z+G))\to \mO_W(-k(K_W+f_*^{-1}\Delta_Z+G))\to \mO_G(-k(K_G+\Delta_G))\to 0. $$
Note that 
$$-G-k(K_W+f_*^{-1}\Delta_Z+G)-(K_W+f_*^{-1}\Delta_Z)=-(k+1)(K_W+f_*^{-1}\Delta_Z+G)$$
is ample over $Z$, by an lc version of Kawamata-Viehweg vanishing theorem (see \cite[Theorem 1.7]{Fujino14}), we have the following surjective map:
$$H^0(W, -k(K_W+f_*^{-1}\Delta_Z+G))\to H^0(G, -k(K_G+\Delta_G))\to 0. $$
We may assume that $k$ is sufficiently large such that $-k(K_G+\Delta_G)$ is very ample, thus one can choose a general element $A_G$ in the linear system $|-k(K_G+\Delta_G)|$ such that $(G, \Delta_G+A_G)$ is still semi-log canonical. Via the above surjective map, one can extend $A_G$ to a divisor $A$ on $W$ such that $A|_G=A_G$. By the inversion of adjunction, the pair $(W, f_*^{-1}\Delta_Z+G+A)$ is log canonical around $G$. Denote by $\Theta:=f_*(\frac{1}{k}A)$, then it is not hard to see the following:
$$K_W+f_*^{-1}\Delta_Z+G+\frac{1}{k}A=f^*(K_Z+\Delta_Z+\Theta). $$
By the choice of $k$, we see that $k\Theta$ is a Cartier divisor on $Z$, then we may write $k\Theta={\rm{div}} (h)$ for some $h\in \mO_{Z,o}$. Note here that $(Z, \Delta_Z+\Theta)$ is log canonical around $o\in Z$ and $v_0$ is an lc place of $(Z, \Delta_Z)$, thus $v_0$ is also an lc place of $(Z, \Delta_Z+\Theta)$. Denote by $\Theta':= \frac{1}{k}\cdot {\rm{div}} (\textbf{in}(h))$, where $\textbf{in}(h)$ is the initial part of $h$, then by \cite[Theorem 3.1]{dFEM10}, the pair $(Z, \Delta_Z+\Theta')$ is also log canonical around $o\in Z$. Since $\Theta'$ is $\bC^*$-equivariant, it can be obtained as the extension of a divisor $D$ on $X$, or in other words, the pair $(Z, \Delta_Z+\Theta')$ is the cone over $(X, \Delta+D)$ with respect to $rL$. By our construction of $\Theta'$, it is not hard to see that $\Theta'$ is an effective divisor which is not zero, and $K_Z+\Delta_Z+\Theta'$ is $\bQ$-Cartier. By \cite[Proposition 3.14(4)]{Kollar13}, we have 
$$L\sim_\bQ a\cdot (K_X+\Delta+D) $$
for some $a\in \bQ^+$. This is a contradiction, since the cone over a stable pair can never be log canonical (see \cite[Lemma 3.1]{Kollar13}). The contradiction implies that $A_{Z, \Delta_Z}(G)=0$, and thus $E$ is an lc place of $(X, \Delta)$.
\end{proof}

\begin{remark}
This remark is pointed out by one of the Referees. In the above proof, by assuming $A_{Z, \Delta_Z}(G)>0$, we obtain a log canonical pair $(Z, \Delta_Z+\Theta)$ such that $v_0$ is an lc place. As $v_0$ is already the lc place of $(Z, \Delta_Z)$ and $\Theta$ passes through the cone vertex via construction, thus $v_0$ cannot be an lc place of $(Z, \Delta_Z+\Theta)$. This already deduces a contradiction which implies that $A_{Z, \Delta_Z}(G)=0$ and we do not actually need the argument on initial degeneration. 
We thank the referee for this simplified argument.
\end{remark}

We next establish the following local correspondence, which is similar to the local correspondence for log Fano cone singularities (see \cite[Proposition 4.25]{Xu21} or \cite[Lemma 2.21]{LWX21}).

\begin{theorem}\label{thm: local correspondence}
There is a  correspondence between non-trivial weakly special test configurations with integral central fibers of $(X, \Delta; L)$ and non-trivial $\bC^*$-equivariant lc blow-ups\footnote{This means that the blow-up is not obtained via canonical blow-up of $o\in Z$ which gives rise to $v_0$.} over $o\in Z$.  More precisely,  given a weakly special test configuration with integral central fiber $(\mX,\Delta_\tc; \mL)\to \bA^1$ of $(X,\Delta; L)$ such that ${\ord_{\mX_0}}|_{K(X)}=a\cdot \ord_E$ for some $a\in \bZ^+$, then one could obtain a $\bC^*$-equivariant lc blowup of the form $\lambda\cdot v_0+a\cdot \ord_{E_\infty}$  for some $\lambda\in \bZ^+$,; conversely, given a non-trivial $\bC^*$-equivariant lc blow-up of the form $\lambda\cdot v_0+c\cdot \ord_{E_\infty}$ for some $\lambda\in \bZ^+$ and $c\in \bQ^+$, then for any given $a\in \bZ^+$ one could construct a weakly special test configuration with integral central fiber $(\mX,\Delta_\tc; \mL)$ such that ${\ord_{\mX_0}}|_{K(X)}=a\cdot \ord_E$.

\end{theorem}

\begin{proof}
Given a non-trivial $\bC^*$-equivariant lc blow-up over $o\in Z$, denoted by $f: W\to Z$, and let $G$ be the exceptional divisor. By the proof of Theorem \ref{thm: key}, there is a prime divisor $E$ over $X$ such that then $\ord_G$ is a  quasi-monomial combination of $v_0$ and $\ord_{E_\infty}$ with weight $(\lambda, c)$ for some $\lambda\in \bZ^{+}$ and $c\in \bQ^{+}$,
where $E_\infty$ is the natural extension of $E$ and $E$ is an lc place of $(X,\Delta)$. The proofs of Lemma \ref{dreamy lemma} and Proposition \ref{lc place to dege} give the way to construct the weakly special test configuration with integral central fiber $(\mX, \Delta_\tc; \mL)$ of $(X, \Delta; L)$ satisfying ${\ord_{\mX_0}}|_{K(X)}=c'\cdot \ord_E$ for any $c'\in \bZ^+$.
It remains to consider the converse direction.

Suppose we are given a weakly special test configuration with integral central fiber $(\mX, \Delta_\tc; \mL)$ of $(X,\Delta;L)$ such that 
$$v_{\mX_0}:={\ord_{\mX_0}}|_{K(X)}=a\cdot \ord_E$$ 
for some prime divisor $E$ over $X$ and $a\in \bZ^+$, we aim to create a $\bC^*$-equivariant lc blow-up associated to $E$ . Denote by
$$R_m:= H^0(X, mrL) \quad \text{and}\quad \mF^j_{v_{\mX_0}}R_m:=\{s\in R_m\mid v_{\mX_0}(s)\geq j\}.$$
By \cite{BHJ17}, we have the following reformulation:
$$\mX=\Proj\ \bigoplus_{m\in \bZ}\bigoplus_{j\in \bZ} \mF^j_{v_{\mX_0}}R_m\cdot t^{-j} \quad \text{and}\quad \mL=\frac{1}{r}\mO_\mX(1).$$
Let us use the following notation for convenience:
$$\mR:= \bigoplus_{m\in \bZ}\bigoplus_{j\in \bZ} \mF^j_{v_{\mX_0}}R_m\cdot t^{-j}.$$
It is clear that $\mR$ is a finitely generated $\bC[t]$-algebra with two gradings, $m$ and $j$. Thus $\mR$ admits a $\bC^*\times \bC^*$-action, where the relative cone structure corresponds to the action by the co-weight $(1,0)$, and $\bC^*$-action from the test configuration corresponds to the co-weight $(0,1)$. Denote by 
$$w_1:=v_0+a\cdot \ord_{E_\infty}, \ \ka_{m,p}:=\{s\in R_m\mid w_1(s)\geq p \},\  \ka_p:=\{s\in R\mid w_1(s)\geq p \} ,$$
then 
$$\ka_p=\bigoplus_{m\in \bN}\ka_{m,p} \ \quad \text{and}\quad \  \ka_{m,p}= \mF^{p-m}_{v_{\mX_0}}R_m.$$
Thus we have the following
$$\bigoplus_{p\in \bZ}\ka_p\cdot t^{-p}=\bigoplus_{p\in \bZ}\bigoplus_{m\in \bZ}\mF^{p-m}_{v_{\mX_0}}R_m\cdot t^{-p}\cong  \bigoplus_{m\in \bZ}\bigoplus_{j\in \bZ} \mF^j_{v_{\mX_0}}R_m\cdot t^{-j}.$$
Therefore,  $\bigoplus_{p\in \bZ}\ka_p$ is finitely generated and the central fiber of the test configuration can be expressed as
$$\mX_0:=\Proj\ \bigoplus_{m\in \bZ}\bigoplus_{j\in \bZ}  \mF^j_{v_{\mX_0}}R_m/\mF^{j+1}_{v_{\mX_0}}R_m\cong \Proj \ \bigoplus_{p\in \bZ}\ka_{p}/\ka_{p+1}.$$
Take a weighted blow-up of $o\in Z$ with respect to the filtration induced by $w_1$, denoted by $\mu_1: W_1\to Z$. More concretely, 
$$W_1=\Proj_Z \bigoplus_{p\in \bZ}\ka_p=\Proj_Z \bigoplus_{m\in \bZ}\bigoplus_{p\in \bZ}\ka_{m,p}.$$
Let $E_1$ be the exceptional divisor of $\mu_1$, then by \cite[Lemma 1.13]{BHJ17},
\begin{align}\label{equ: pl}
E_1=\Proj\ \bigoplus_{p\in \bZ}\ka_{pl}/\ka_{pl+1}
\end{align}
for a sufficiently divisible positive integer $l$, and $-E_1$ is ample over $Z$. 
Consider the following degeneration of $Z$:
$$\mZ:=\Spec\ \bigoplus_{p\in \bZ}\ka_p\cdot t^{-p}\to \bA^1 $$
with the central fiber
\begin{align}\label{equ: p}
\mZ_0:=\Spec\ \bigoplus_{p\in \bZ}\ka_p/\ka_{p+1} ,
\end{align}
and denote by $\Delta_\mZ$ the natural extension of $\Delta_Z$, then $(\mZ, \Delta_\mZ+\mZ_0)$ is log canonical since $(\mX, \Delta_\tc+\mX_0)$ is log canonical. By sub-adjunction we see that $(\mZ_0, \Delta_{\mZ_0})$ is slc, where $\Delta_{\mZ_0}:={\Delta_\mZ}|_{\mZ_0}$. Write 
$$K_{E_1}+\Delta_{E_1}:=(K_W+{\mu_1}_*^{-1}\Delta_Z+E_1)|_{E_1}, $$
then by (\ref{equ: pl}) and (\ref{equ: p}) we see that  $(\mZ_0,\Delta_{\mZ_0})$ is the orbifold cone over $(E_1, \Delta_{E_1})$ and thus $(E_1, \Delta_{E_1})$ is slc.
By the inversion of adjunction,  $\mu_1: (W_1, E_1)\to Z $
is a $\bC^*$-equivariant lc blow-up and the proof is finished. It is worth to say one more sentence. Recall that $\ord_{E_1}$ is a quasi-monomial combination of $v_0$ and $\ord_{E_\infty}$ with weight $(1, a)$ and $v_{\mX_0}=a\cdot \ord_E$, by Theorem \ref{thm: key}, $E$ is an lc place of $(X,\Delta)$.
\end{proof}

\

\section{Global correspondence}\label{sec: global correspondence}

In this section, we prove Theorem \ref{thm: main1}. We begin with the following dreamy property.

\begin{lemma}\label{dreamy lemma}
Let $(X,\Delta)$ be a projective lc CY pair and $L$ an ample $\bQ$-line bundle on $X$. If $E$ is an lc place of $(X,\Delta)$, then $E$ is dreamy with respect to $L$.
\end{lemma}

\begin{proof}
Choose a divisible positive integer $r$ such that $rL$ is Cartier. Let $Z$ be the affine cone over $X$ with respect to $rL$, i.e.,
$$Z=\Spec \bigoplus_{m\in \bN}H^0(X, mrL),$$
and denote by $o\in Z$ the cone vertex. Let $\Delta_Z$ be the extension of $\Delta$ on $Z$. Consider the projection morphism $p: Z\setminus o\to X$, we denote $E_\infty$ to be the prime divisor over $Z\setminus o$ via pulling back $E$ over $X$. Let $\mu: \tZ\to Z$ be the blow-up of the vertex and $X_0$ the exceptional divisor. Write $v_0:=\ord_{X_0}$ for the canonical valuation. It is clear that $E_\infty$ is an lc place of the pair $(Z, \Delta_Z)$, and by \cite[Proposition 3.14]{Kollar13}, $v_0$ is also an lc place of $(Z, \Delta_Z)$. 

Let $w_k:=k\cdot v_0+\ord_{E_\infty}$ be the quasi-monomial valuation with weight $(k,1)$ along $v_0$ and $\ord_{E_\infty}$ for some positive integer $k$. By \cite[Proposition 5.1]{JM12}, $w_k$ is an lc place of $(Z,\Delta_Z)$ whose center is exactly the cone vertex $o\in Z$. By \cite[Lemma 5.1]{CZ21}, there is an extraction morphism $f: W\to Z$ which has a unique exceptional divisor corresponding to $w_k$, denoted by $E_k\subset W$, such that $-E_k$ is ample over $Z$, and we have the following:
$$K_W+f_*^{-1}\Delta_Z+E_k=f^*(K_Z+\Delta_Z). $$
Consider the following exact sequence for $pl\in \bN$, where $l\in \bN$ and $p$ is a fixed divisible positive integer such that $pE_k$ is Cartier:
$$0\to \mO_{W}(-(pl+1)E_k)\to \mO_{W}(-plE_k)\to \mO_{E_k}(-plE_k)\to 0. $$
Write $\ka_l:=f_*\mO_{W}(-lE_k)$, then 
$$\ka_{pl}/\ka_{pl+1}\cong H^0(E_k, -plE_k|_{E_k})$$ 
for $l\in \bN$, since $-E_k$ is $f$-ample and $R^1 f_*\mO_{W}(-(pl+1)E_k)=0$ by an lc version of Kawamata-Viehweg vanishing (see \cite[Theorem 1.7]{Fujino14}). Thus the graded algebra $\bigoplus_{l\in \bN}\ka_{pl}/\ka_{pl+1}$ is finitely generated, hence so is $\bigoplus_{l\in \bN}\ka_l/\ka_{l+1}$.  Therefore, the graded algebra $\bigoplus_{l\in \bN}\ka_l$ is finitely generated.

Denote by 
$$R_m:=H^0(X,mrL) \quad \text{and}\quad  R:=\bigoplus_{m\in \bN}R_m.$$ 
We naturally extend $\bigoplus_{l\in \bN}\ka_l$ and $R$ to graded algebras indexed by $\bZ$ via defining $R_m=0$ for $m<0$ and $\ka_l=\mO_Z$ for $l< 0$. We use the notation 
$$\mF^j_{\ord_E}R_m:=\{s\in R_m\mid  \textit{$\ord_E(s)\geq j$}\}\quad \text{and} \quad \ka_{m,l}:=\{s\in R_m\mid \textit{$km+\ord_E(s)\geq l$}\} .$$
It is clear that $\ka_l=\bigoplus_{m\in \bZ}\ka_{m,l}$, and $s\in \mF^j_{\ord_E}R_m$ if and only if $s\in \ka_{m,km+j}$.
Then we have
$$\bigoplus_{m\in \bZ}\bigoplus_{j\in \bZ} \mF^j_{\ord_E}R_m\cong \bigoplus_{m\in \bZ}\bigoplus_{j\in \bZ} \ka_{m,km+j}\cong \bigoplus_{m\in \bZ}\bigoplus_{j\in \bZ} \ka_{m,j}\cong \bigoplus_{j\in \bZ} \ka_j,$$
which are all finitely generated graded algebras.
Note that 
$$\mF^j_{\ord_E}R_m\cong H^0(Y, g^*(mrL)-jE ), $$
where $g: Y\to X$ is a projective normal birational model such that $E$ is a prime divisor on $Y$ (see \cite[Lemma 5.1]{CZ21}). Then $\bigoplus_{m\in \bN}\bigoplus_{j\in \bN} \mF^j_{\ord_E}R_m$ being finitely generated implies that $E$ is dreamy with respect to $L$. The proof is finished.
\end{proof}

\begin{remark}
The idea of the above proof is the same as that of \cite[Theorem 1.5]{CZ21}, and \cite[Theorem 1.5]{CZ21} is just a direct corollary of Lemma \ref{dreamy lemma}.
\end{remark}

Next we show that an lc place induces a weakly special degeneration.

\begin{proposition}\label{lc place to dege}
Let $(X,\Delta)$ be a projective lc CY pair and $L$ an ample $\bQ$-line bundle on $X$. Suppose $E$ is an lc place of $(X,\Delta)$, then for any given $c\in \bZ^+$, there exists a non-trivial weakly special test configuration with integral central fiber $(\mX, \Delta_\tc; \mL)$ of $(X,\Delta; L)$ such that $\ord_{\mX_0}|_{K(X)}=c\cdot \ord_E$.
\end{proposition}

\begin{proof}
By \cite[Lemma 5.1]{CZ21}, there exists an extraction $g: Y\to X$ which only extracts $E$, and we have
$$K_Y+g_*^{-1}\Delta+E=g^*(K_X+\Delta). $$
If $E$ is a prime divisor on $X$, then we may just take $g={\rm id}$. 
We only deal with the case where $E$ is exceptional in the remaining proof, as one can prove similarly for the case where $E$ is a divisor on $X$.
According to Lemma \ref{dreamy lemma}, we know that the following $\bN^2$-graded ring 
$$\bigoplus_{k\in \bN}\bigoplus_{j\in \bN} H^0(Y, g^*(krL)-jE) $$
is finitely generated for some $r\in \bZ^+$ such that $rL$ is Cartier. Let 
$$R_k:=H^0(X, krL) \ \text{ and }\  \mF^j_{c\cdot \ord_E}R_k:=\{s\in R_k\mid c\cdot \ord_E(s)\geq j\}.$$ 
We construct the following degeneration family over $\bA^1$:
$$\mX:=\Proj \bigoplus_{k\in \bN}\bigoplus_{j\in \bZ} \left(\mF^j_{c\cdot \ord_E}R_k\right) t^{-j}\to \bA^1.$$
 By the same proof as that of \cite[Lemma 3.8]{Fuj17}, we know that $(\mX, \Delta_\mX; \mO_\mX(1))$ is a test configuration of $(X,\Delta; rL)$ whose central fiber $\mX_0$ is integral, where  $\Delta_\mX$ is the extension of $\Delta$  on $\mX$.  Then it holds that ${\ord_{\mX_0}}|_{K(X)}= c\cdot \ord_E$.
We next show that $(\mX, \Delta_{\mX}+\mX_0)$ is log canonical, which will imply that $(\mX, \Delta_{\mX}; \frac{1}{r}\mO_\mX(1))$ is a weakly special test configuration of $(X,\Delta; L)$ induced by $E$.

Put the morphism $g$ into the trivial family $G: Y_{\bA^1}\to X_{\bA^1}$, then we have
$$K_{Y_{\bA^1}}+G_*^{-1}\Delta_{\bA^1}+E_{\bA^1}+Y_0=G^*\left(K_{X_{\bA^1}}+\Delta_{\bA^1}+X_0\right) ,$$
where $X_0$ (resp. $Y_0$) is the central fiber of the family $X_{\bA^1}\to \bA^1$ (resp. $Y_{\bA^1}\to \bA^1$). 
Let $v$ be the divisorial valuation over $X_{\bA^1}$ with weight $(c,1)$ along divisors $E_{\bA^1}$ and $Y_0$, then 
$$A_{X_{\bA^1}, \Delta_{\bA^1}+X_0}(v)=A_{Y_{\bA^1}, G_*^{-1}\Delta_{\bA^1}+E_{\bA^1}+Y_0}(v)=0. $$
By \cite[Lemma 5.1]{CZ21}, one can extract $v$ to a divisor on a birational model $h: \mY\to X_{\bA^1}$ as follows:
$$K_{\mY}+h_*^{-1}\Delta_{\bA^1}+h_*^{-1}X_0+\mE=h^*(K_{X_{\bA^1}}+\Delta_{\bA^1}+X_0),$$
where $v=\ord_\mE$.
We note here that the pair $\left(\mY, h_*^{-1}\Delta_{\bA^1}+h_*^{-1}X_0+\mE\right)$ is  log canonical and log Calabi-Yau over $\bA^1$. 

Now we could compare the two pairs, 
$$\left(\mY, h_*^{-1}\Delta_{\bA^1}+h_*^{-1}X_0+\mE\right) \quad \text{and}\quad (\mX, \Delta_{\mX}+\mX_0).$$ 
Since $\mE$ and $\mX_0$ induce the same divisorial valuation on $K(\mY)\cong K(X\times \bA^1)\cong K(\mX)$, we have the following birational contraction map:
$$\left(\mY, h_*^{-1}\Delta_{\bA^1}+h_*^{-1}X_0+\mE\right)\dashrightarrow (\mX, \Delta_{\mX}+\mX_0). $$
As both two sides are log Calabi-Yau over $\bA^1$, they are crepant, i.e., for any common log resolution $p: \mW\to \mY$ and $q: \mW\to \mX$ one has
$$p^*\left(K_{\mY}+h_*^{-1}\Delta_{\bA^1}+h_*^{-1}X_0+\mE\right)=q^*(K_{\mX}+\Delta_{\mX}+\mX_0). $$
It follows that $(\mX, \Delta_{\mX}+\mX_0)$ is log canonical, and this completes the proof.
\end{proof}

We turn to the converse direction.

\begin{proposition}\label{dege to lc place}
Let $(X,\Delta)$ be a projective lc CY pair and $L$ an ample $\bQ$-line bundle on $X$. Given a non-trivial weakly special test configuration with integral central fiber $(\mX, \Delta_\tc; \mL)$ of $(X,\Delta; L)$ such that ${\ord_{\mX_0}}|_{K(X)}=c\cdot \ord_E$ for some prime divisor $E$ over $X$ and $c\in \bZ^+$,  then $E$ is an lc place of $(X,\Delta)$.
\end{proposition}

\begin{proof}
By the proof of Theorem \ref{thm: local correspondence}, there is a $\bC^*$-equivariant lc blow-up over $o\in Z$, denoted by $\mu_1: W_1\to Z$, such that
$$\ord_{E_1}=c\cdot \ord_{E_\infty}+v_0, $$
where $E_1$ is the exceptional divisor of $\mu_1$, and $E_\infty$ is the natural extension of $E$.  By Theorem \ref{thm: key} we see that $E$ is an lc place of $(X,\Delta)$.
\end{proof}

\begin{proof}[Proof of Theorem \ref{thm: main1}]
The proof is just a combination of Lemma \ref{dreamy lemma}, Proposition \ref{lc place to dege} and Proposition \ref{dege to lc place}.
\end{proof}

\


\section{Generalized global correspondence}\label{sec: g-global correspondence}

In this section, we will generalize Theorem \ref{thm: main1} for weakly special test configurations without assuming the irreducibility of central fibers. 



Let $(X, \Delta)$ be a projective lc CY pair and $L$ an ample $\bQ$-line bundle on $X$. Given a weakly special semi-test configuration $(\mX, \Delta_\tc; \mL)$ of $(X, \Delta; L)$, and denote by 
$$v_{\mX_{0,i}}:={\ord_{\mX_{0,i}}}|_{K(X)}=c_i\cdot \ord_{E_i}, $$
where $\mX_{0,i}, i=1,...,k,$ are all non-trivial components\footnote{Here we say that a component in $\mX_0$ is non-trivial if it is not the strict transform of $X\times 0$ with respect to the birational map $\mX\dashrightarrow X\times \bA^1$, see \cite[Definition 4.4]{BHJ17}. } in the central fiber $\mX_0$, $c_i\in \bZ^+$, and $E_i$ are prime divisors over $X$. We have the following characterization of $E_i$.

\begin{theorem}\label{thm: g-key}
Notation as above. Then $E_i$ is an lc place of $(X, \Delta)$ for every $i$.
\end{theorem}

\begin{proof}
Let 
$$h: (\mY, \Delta_\mY+\mY_0)\to (\mX, \Delta_\tc+\mX_0)$$ 
be a $\bC^*$-equivariant $\bQ$-factorial dlt modification such that
$$K_\mY+\Delta_\mY+\mY_0=h^*(K_\mX+\Delta_\tc+\mX_0).$$
Let  $(\mY_t, \Delta_{\mY_t})$ be the fiber over $0\ne t\in \bA^1$, then we see that 
$$\phi: (\mY_t, \Delta_{\mY_t}) \to (X, \Delta)$$
is a crepant lc model of $(X,\Delta)$, and $(\mY, \Delta_\mY; h^*\mL)$ is a semi-test configuration of $(\mY_t, \Delta_{\mY_t}; L' )$, where 
$L'=\phi^*L$. In particular, $(\mY_t, \Delta_{\mY_t})$ is an lc CY pair. Denote $\mY_{0,i}$ to be the strict transform of $\mX_{0,i}$ and we still have ${\ord_{\mY_{0,i}}}|_{K(\mY_t)}=c_i\cdot \ord_{E_i}$, where $E_i$ can also be viewed as a prime divisor over $\mY_t$. It suffices to show that $E_i$ is an lc place of $(\mY_t, \Delta_{\mY_t})$. 

Consider the pair $(\mY, \Delta_\mY+\mY_{0,i})$, it is clear to see the following:
$$K_{\mY}+\Delta_\mY+\mY_{0,i}\sim_{\bQ,\bA^1} -(\mY_0-\mY_{0,i}) .$$
By \cite[Theorem 1.1]{Birkar12} or \cite[Theorem 1.6]{HX13}, one may run a $\bC^*$-equivariant MMP on $K_{\mY}+\Delta_\mY+\mY_{0,i}$ over $\bA^1$, which terminates with a minimal model, denoted by 
$$(\mY, \Delta_\mY+\mY_{0,i})\dashrightarrow (\mY', \Delta_{\mY'}+\mY'_{0,i}), $$
where $\Delta_{\mY'}$ (resp. $\mY'_{0,i}$) is the push-forward of $\Delta_\mY$ (resp. $\mY_{0,i}$). Thus $-(\mY'_0-\mY'_{0,i})$ is nef over $\bA^1$, or equivalently, $\mY'_{0,i}$ is nef over $\bA^1$. By the Zariski lemma (e.g. \cite[Section 2.4]{LX14}), $\mY'_{0,i}=\mY'_0$. It is not hard to see that the following two pairs
$$(\mY, \Delta_\mY+\mY_0)\ \quad \text{and}\quad\  (\mY', \Delta_{\mY'}+\mY'_{0,i})$$
are isomorphic over $\bA^1\setminus 0$ and crepant to each other,  thus $(\mY',\Delta_{\mY'})$ is a weakly special degeneration of $(\mY_t, \Delta_{\mY_t})$ in the sense of Definition \ref{defn: tc}. Put $A:= {\mO_{\mY'}(1)}|_{\mY_t}$, then $(\mY',\Delta_{\mY'}; \mO_{\mY'}(1))$ is a weakly special test configuration with integral central fiber of $(\mY_t, \Delta_{\mY_t}; A)$.
Note that $\ord_{\mY_{0,i}}$ and $\ord_{\mY'_{0}}$ are the same divisorial valuation over the function field 
$$K(\mX)=K(\mY)=K(\mY')=K(X)(t),$$ 
thus $E_i$ is an lc place of $(\mY_t, \Delta_{\mY_t})$ by Theorem \ref{thm: key} and Theorem \ref{thm: local correspondence}. The proof is finished.
\end{proof}

We still denote $(Z, \Delta_Z)$ to be the affine cone over the lc CY pair $(X, \Delta)$ with respect to $rL$, and $o\in Z$ is the cone vertex. 

\begin{definition}\label{g-lc blowup}
We say that a finite set of divisorial valuations over $X$, denoted by 
$$\{c_1\cdot \ord_{E_1},...,c_l\cdot \ord_{E_l}\},$$ 
is a \emph{weakly special collection} of $(X, \Delta)$ if $c_i\in \bZ^+$ and $E_i$ is an lc place of $(X, \Delta)$ for every $i$.
\end{definition}

Note here that our definition of weakly special collection is a little different from that in \cite[Appendix A]{BLX19}, since all $w_i$ are automatically induced by lc places of $(X,\Delta)$. However, we will later see that these definitions are essentially equivalent. We aim to establish the following generalized global correspondence.

\begin{theorem}\label{thm: g-global correspondence}
Let $(X, \Delta)$ be a projective lc CY pair and $L$ an ample $\bQ$-line bundle on $X$, then there exists a  correspondence between weakly special collections of $(X, \Delta)$ and weakly special semi-test configurations of $(X, \Delta; L)$. More precisely, given a weakly special collection $\{c_1\cdot \ord_{E_1},...,c_l\cdot \ord_{E_l}\},$ we can construct a weakly special semi-test configuration  
$(\mX, \Delta_\tc; \mL)$ satisfying 
$${\ord_{\mX_{0,i}}}|_{K(X)}=c_i\cdot \ord_{E_i}, i=1,...,l,$$
where $\mX_{0,i}$ are all non-trivial components of $\mX_0$; conversely, given a weakly special semi-test configuration  
$(\mX, \Delta_\tc; \mL)$ satisfying 
$${\ord_{\mX_{0,i}}}|_{K(X)}=c_i\cdot \ord_{E_i}, i=1,...,l,$$
where $c_i\in \bZ^+$ and $\mX_{0,i}$ are all non-trivial components of $\mX_0$, then the set $\{c_1\cdot \ord_{E_1},...,c_l\cdot \ord_{E_l}\}$ is a weakly special collection.
\end{theorem}

\begin{proof}
Given a weakly special semi-test configuration of the form in the theorem, denoted by $(\mX, \Delta_\tc; \mL)$, by Theorem \ref{thm: g-key}, every $E_i$ is an lc place of $(X, \Delta)$. Thus $\{c_i\cdot \ord_{E_i}\}_i$ is a weakly special collection of $(X, \Delta)$ induced by the given weakly special test configuration.

Now we consider the converse direction. Given a weakly special collection 
$$\{c_1\cdot \ord_{E_1},...,c_l\cdot \ord_{E_l}\},$$ 
we aim to construct a weakly special semi-test configuration as in the theorem. Let $G_i(\lambda_i)$ be the prime divisor over $Z$ such that 
$$\ord_{G_i(\lambda_i)}:=\lambda_i\cdot v_0+c_i\cdot \ord_{E_{i,\infty}},$$ 
where $\lambda_i\in \bZ^+$. Then all $G_i(\lambda_i)$ are lc places of $(Z, \Delta_Z)$ centered at $o\in Z$. We have the following claim:

\begin{claim}\label{claim: extraction}
There exists a birational model $\mu: W\to Z$ which  extracts part of $G_i(\lambda_i)$, say $\{G_{j_1}(\lambda_{j_1}),..., G_{j_{l'}}(\lambda_{j_{l'}})\}$, such that $-\sum_{s=1}^{l'} G_{j_s}(\lambda_{j_s})$ is $\mu$-ample. Here $\{j_1,...,j_{l'}\}$ is a non-repeating subset of $\{1,...,l\}$.
\end{claim}
\begin{proof}[Proof of the Claim]
Let $\gamma: Y\to Z$ be a $\bC^*$-equivariant $\bQ$-factorial dlt modification such that all $G_i(\lambda_i)$ appear as prime divisors on $Y$. We write
$$K_Y+\gamma_*^{-1}\Delta_Z+\sum_{i=1}^{l} G_i(\lambda_i)+C=\gamma^*(K_Z+\Delta_Z), $$
where $C$ is the reduced sum of lc places other than $\sum_{i=1}^lG_i(\lambda_i)$. Choose a small $0<\epsilon \ll 1$ and consider the pair 
$$(Y, \gamma_*^{-1}\Delta_Z+(1-\epsilon)\cdot \sum_{i=1}^{l} G_i(\lambda_i)+C),$$
then we have
$$K_Y+\gamma_*^{-1}\Delta_Z+(1-\epsilon)\cdot \sum_{i=1}^{l} G_i(\lambda_i)+C\sim_{\gamma,\bQ} -\epsilon \sum_{i=1}^{l} G_i(\lambda_i) .$$
By applying \cite[Theorem 1.1]{Birkar12} or \cite[Theorem 1.6]{HX13} and arguing by the same way as \cite[Lemma 5.1]{CZ21}, one can run a $\bC^*$-equivariant MMP over $Z$ on 
$$K_Y+\gamma_*^{-1}\Delta_Z+(1-\epsilon)\cdot \sum_{i=1}^{l} G_i(\lambda_i)+C$$ 
to get a good minimal model over $Z$, denoted by
$$Y\dashrightarrow W'/Z,$$
and the ample model $\mu: W\to Z$ is exactly what we want.
Consider the birational contraction $Y\dashrightarrow W$, by the same argument of \cite[Lemma 5.1]{CZ21}, we know that 
the components in $\sum_{i=1}^{l} G_i(\lambda_i)$ cannot be all  contracted (but some may be contracted),  and this is why we say that only part of $G_i(\lambda_i)$ are extracted in the statement of the Claim.
For convenience, we just assume $\{j_1,...,j_{l'}\}=\{1,...,l'\}$ for some $l'<l$, and still use $\sum_{i=1}^{l'} G_{i}(\lambda_{i})$ to denote the push-forward of $\sum_{i=1}^{l} G_{i}(\lambda_{i})$ via $Y\dashrightarrow W$. Thus $-\sum_{i=1}^{l'} G_{i}(\lambda_{i})$ is $\mu$-ample.\end{proof}

We first construct a weakly special test configuration  
$(\mX, \Delta_\tc; \mL)$ satisfying 
$${\ord_{\mX_{0,i}}}|_{K(X)}=c_i\cdot \ord_{E_i}, i=1,...,l',$$
where $\mX_{0,i}, i=1,...,l',$ are all non-trivial components of $\mX_0$.
Recall that 
$$Z=\Spec\ R=\Spec\ \bigoplus_{m\in \bN}R_m=\Spec\ \bigoplus_{m\in \bN} H^0(X, mrL),$$ 
we first introduce the following ideals for $i=1,...,l'$:
$$\ka_{m,p}^{(i)}:=\{s\in R_m\mid \ord_{G_i(\lambda_i)}(s)=c_i\cdot \ord_{E_i}(s)+\lambda_i\cdot v_0(s)\geq p\}, $$
$$\ka_p^{(i)}:=\{s\in R\mid \ord_{G_i(\lambda_i)}(s)=c_i\cdot \ord_{E_i}(s)+\lambda_i\cdot v_0(s)\geq p\},$$
$$I_{m,p}:=\cap_i \ \ka_{m,p}^{(i)} \quad \text{and} \quad I_p:=\cap_i \ \ka_p^{(i)}.$$
Note that $-\sum_{i=1}^{l'} G_i(\lambda_i)$ is ample over $Z$ (see  Claim \ref{claim: extraction}), 
thus by \cite[Lemma 1.8  and Theorem 1.10]{BHJ17}, we have the following formulation of $W$: 
$$W=\Proj_Z\ \bigoplus_{p\in \bZ}I_p=\Proj_Z \bigoplus_{m\in \bN}\bigoplus_{p\in \bZ}I_{m,p}. $$
Note that the ideal sequence $I_\bullet$ induces a degeneration of $Z$ over $\bA^1$ via
$$\mZ:=\Spec\ \bigoplus_{p\in \bZ} I_p\cdot t^{-p}=\Spec\ \bigoplus_{m\in \bN}\bigoplus_{p\in \bZ} I_{m,p}\cdot t^{-p}\to \bA^1, $$
where the central fiber is given by
$$\mZ_0:=\Spec\ \bigoplus_{p\in \bN}I_p/I_{p+1} =\Spec\ \bigoplus_{m\in \bN}\bigoplus_{p\in \bN}I_{m,p}/I_{m,p+1}.$$
Let $\Delta_{\mZ_0}$ be the restriction of $\Delta_\mZ$ on $\mZ_0$, where $\Delta_{\mZ}$ is the extension of $\Delta_Z$ on $\mZ$. Then we know that $(\mZ, \Delta_\mZ,\xi;\eta)$ is a test configuration of $(Z,\Delta_Z, \xi)$ in the sense of \cite[Definition 2.14]{LWX21}, where $\xi$ (resp. $\eta$) is the vector field on $Z$ (resp. $\mZ$) induced by the grading $m$ (resp. $p$).
We have the following claim:
\begin{claim}\label{claim: cone}
Notation as above, $(\mZ, \Delta_\mZ,\xi;\eta)$ is a weakly special test configuration of $(Z,\Delta_Z, \xi)$ in the sense of \cite[Definition 2.14]{LWX21}, more precisely, $(\mZ_0, \Delta_{\mZ_0})$ is semi-log canonical and it is an orbifold cone over an slc CY pair.
\end{claim}
\begin{proof}[Proof of the Claim]
Denote by $G:=\sum_{i=1}^{l'} G_i(\lambda_i)$, then $-G$ is ample over $Z$ as we have mentioned. Choose a divisible $k\in \bN$ such that $-kG$ is Cartier. By \cite[Lemma 1.8 and Theorem 1.10]{BHJ17}, we have the following characterization for any $p\in \bZ$:
$$I_{p}=\mu_*\mO_W(-pG).$$
Consider the $\mu_k$-action on $(\mZ, \Delta_\mZ)$, where $\mu_k$ is the multiplicative group of $k$-th roots of unity. Let $(\mZ', \Delta_{\mZ'}):=(\mZ, \Delta_\mZ)/\mu_k$, then we have
$$\mZ'=\Spec\ \bigoplus_{p\in \bZ}I_{kp}\cdot t^{-p}\to \bA^1,$$
and the quotient map $\sigma: \mZ\to \mZ'$ is a lifting of the map $\bA^1\to \bA^1, t\mapsto t^k$. It is clear that $\sigma$ is \'etale away from the central fiber. Recall that 
$$\mZ_0=\Spec\ \bigoplus_{p\in \bN}I_p/I_{p+1},$$
thus
$$\mZ_0/\mu_k=\Spec\ \bigoplus_{p\in \bN}I_{kp}/I_{kp+1} \quad \text{and} \quad \Supp\ \mZ_0/\mu_k=\Supp \ \mZ_0'.$$
Let $\red(\mZ_0)$ (resp. $\red(\mZ_0')$) be the reduced part of $\mZ_0$ (resp. $\mZ_0'$), then clearly we have
$$K_\mZ+\Delta_\mZ+\red(\mZ_0) = \sigma^*(K_{\mZ'}+\Delta_{\mZ'}+\red(\mZ_0')).$$
Hence it suffices to show the following three points to complete the proof:
\begin{enumerate}
\item the pair $(\mZ', \Delta_{\mZ'}+\red(\mZ_0'))$ is log canonical;
\item $\mZ_0=\red(\mZ_0)$, i.e., $\mZ_0$ is reduced;
\item $(K_{\mZ'}+\Delta_{\mZ'}+\red(\mZ_0'))|_{\red(\mZ_0')}$ is an affine cone over an slc CY pair.
\end{enumerate}
The first two points together imply that $(\mZ, \Delta_\mZ+\mZ_0)$ is log canonical (e.g. \cite[Proposition 5.20]{KM98}), and the third point implies that $(\mZ_0, \Delta_{\mZ_0})$ is an orbifold cone over some slc CY pair.

We first show that $\mZ_0/\mu_k=C_a(G, \mO_G(-kG))$, where $C_a(G, \mO_G(-kG))$ is the affine cone over $G$ with respect to $-kG|_G$. Consider the following exact sequence
$$0\to \mO_W(-(kp+1)G)\to \mO_W(-kpG)\to \mO_G(-kpG)\to 0. $$
By Claim \ref{claim: extraction} we have
$$K_W+\mu_*^{-1}\Delta_Z+G=\mu^*(K_Z+\Delta_Z), $$
where both two sides are lc CY pairs. Thus  by an lc version of Kawamata-Viehweg vanishing theorem (see \cite[Theorem 1.7]{Fujino14}), we see
$$R^1\mu_*{\mO_W(-(kp+1)G)}=0$$ 
for any $p\in \bN$. Therefore, we have 
$$I_{kp}/I_{kp+1}\cong H^0(G, -kpG). $$
and hence 
$$\red(\mZ_0')=\mZ_0/\mu_k= \Spec\ \bigoplus_{p\in \bN}I_{kp}/I_{kp+1}=\Spec\ \bigoplus_{p\in \bN}\ H^0(G, -kpG).$$
Denote by 
$$K_G+\Delta_G:= (K_W+\mu_*^{-1}\Delta_Z+G)|_G,$$ 
which by adjunction is an slc CY pair, thus the affine cone $(\red(\mZ_0'), \Delta_{\red(\mZ_0')})$ over $(G, \Delta_G)$ with respect to $-kG|_G$ is slc, where 
$$K_{\red(\mZ_0')}+\Delta_{\red(\mZ_0')}:=(K_{\mZ'}+\Delta_{\mZ'}+\red(\mZ_0'))|_{\red(\mZ_0')}.$$
This proves the third point. Apply the inversion of adjunction we conclude the first point. It remains to show the second point. Recall the following characterization:
\begin{align*}
\bigoplus_{p\in \bN}I_p/I_{p+1}= \bigoplus_{p\in \bN} \mu_*{\mO_W(-pG)}/\mu_*{\mO_W(-(p+1)G)}.
\end{align*}
Suppose $0\ne f\in I_p/I_{p+1}$, then there must exist some $i$ such that $\ord_{G_i}(f)=p$, thus $\ord_{G_i}(f^b)=bp$ for any $b\in \bZ^+$. This implies that $0\ne f^b\in I_{bp}/I_{bp+1}$ for any $b\in \bZ^+$ and hence $\mZ_0$ is reduced. The proof of the Claim is finished.
\end{proof}

Define the following degeneration over $\bA^1$:
$$\mX:=\Proj\ \bigoplus_{p\in \bZ} I_p\cdot t^{-p}=\Proj\ \bigoplus_{m\in \bN}\bigoplus_{p\in \bZ} I_{m,p}\cdot t^{-p}\to \bA^1.$$
Then the central fiber can be formulated as follows:
$$\mX_0=\Proj\ \bigoplus_{p\in \bN}I_p/I_{p+1} =\Proj\ \bigoplus_{m\in \bN}\bigoplus_{p\in \bN}I_{m,p}/I_{m,p+1}.$$
It is clear that $(\mX, \mO_\mX(1))$ is a test configuration of $(X, rL)$, where $\mO_\mX(1)$ is the relatively ample line bundle on $\mX$ under the grading $m$. Denote $\Delta_\tc$ the extension of $\Delta$ on $\mX$, then we aim to show that $(\mX, \Delta_\tc; \frac{1}{r}\mO_\mX(1))$ is a weakly special test configuration of $(X,\Delta; L)$. By our construction, $(\mZ_0, \Delta_{\mZ_0})$ is the affine cone over $(\mX_0, \Delta_{\tc,0})$ with respect to $\mO_{\mX_0}(1)$, where $\mO_{\mX_0}(1):={\mO_{\mX}(1)}|_{\mX_0}$. By Claim \ref{claim: cone} we see that $(\mX_0, \Delta_{\tc, 0})$ is semi-log canonical, and thus $(\mX, \Delta_\tc+\mX_0)$ is log canonical by inversion of adjunction. It is not hard to see from the construction that the components in the central fiber exactly produce the divisorial valuations $\{c_1\cdot \ord_{E_1},...,c_{l'}\cdot \ord_{E_{l'}}\}$, which is a part of  the given weakly special collection. 

The final step is to construct a weakly special semi-test configuration  
$(\mX', \Delta'_\tc; \mL')$ satisfying 
$${\ord_{\mX'_{0,i}}}|_{K(X)}=c_i\cdot \ord_{E_i}, i=1,...,l,$$
where $\mX'_{0,i}, i=1,...,l,$ are all non-trivial components of $\mX'_0$.

Take a dlt modification $Y\to X$ which exactly extracts all $E_i, i=1,...l,$ (by the similar argument as \cite[Lemma 5.1]{CZ21} or Claim \ref{claim: extraction}). Then we consider the product space 
$$Y\times \bA^1\to X\times \bA^1.$$
Let us write $X_0, Y_0$  for the central fibers and still denote $E_{i,\infty}$ to be the extension of $E_i$ on $Y\times \bA^1$. Then we get a collection of divisorial valuations
$$\big\{ \lambda_j\cdot\ord_{Y_0}+ c_j\cdot \ord_{E_{j,\infty}}\mid j=1,...,l    \big\}$$
over $X\times \bA^1$. Let $\mE_j$ be a prime divisor over $X\times \bA^1$ such that
$$\ord_{\mE_j}:=\lambda_j\cdot\ord_{Y_0}+ c_j\cdot \ord_{E_{j,\infty}},$$
we see that all $\mE_j$ are log canonical places of $(X\times \bA^1, \Delta\times \bA^1+X_0)$ centered in $X_0$. It is not hard to see that 
$(X\times \bA^1, \Delta\times \bA^1+X_0)$ is crepant to the weakly special test configuration $(\mX, \Delta_{\tc}+\mX_0)$ which we have constructed above, thus $\{\ord_{\mE_j}| l'<j\leq l\}$ are all lc places of $(\mX, \Delta_{\tc}+\mX_0)$.
Apply the similar argument as Claim \ref{claim: extraction}, one can find a 
$\bC^*$-equivariant extraction $\sigma: \mX'\to \mX$ which exactly extracts $\{\ord_{\mE_j}| l'<j\leq l\}$. Let $\Delta'_{\tc}$ be the extension of $\Delta$ on $\mX'$, then $(\mX', \Delta'_{\tc}; \sigma^*\mL)$ is the weakly special semi-test configuration as we want.
\end{proof}

\
\

\begin{remark}\label{rem: set of tc}
In the above proof, given a weakly special collection $\{c_1\cdot \ord_{E_1},...,c_l\cdot \ord_{E_l}\}$ and a vector of positive integers $\vec{\lambda}:=(\lambda_1,...,\lambda_l)\in (\bZ^+)^{\oplus l}$, we define the corresponding prime divisors 
$$\big\{G_1(\lambda_1) ,...,G_l(\lambda_l)\big\}$$ 
over $o\in Z$. Then, we use the sequence of ideals associated to 
$$G:=\sum_{i=1}^{l'} G_i(\lambda_i)$$
to construct the following weakly special test configuration
$$\mX(\vec{\lambda}):=\Proj\ \bigoplus_{p\in \bZ} I_p\cdot t^{-p}=\Proj\ \bigoplus_{m\in \bN}\bigoplus_{p\in \bZ} I_{m,p}\cdot t^{-p}\to \bA^1 $$
with the polarization $\mO_{\mX(\vec{\lambda})}(1)$, which is the relatively ample line bundle on $\mX(\vec{\lambda})$ under the grading $m$. Finally, we construct a weakly special semi-test configuration via a $\bC^*$-equivariant extraction
$$\sigma: \mX(\vec{\lambda})'\to \mX(\vec{\lambda}), $$
which extracts the rest components in the central fiber associated to the given weakly special collection.
This means that we could construct various weakly special semi-test configurations of $(X,\Delta;L)$, 
denoted by 
$$\big\{(\mX(\vec{\lambda})', \Delta_{\mX(\vec{\lambda})'}; \frac{1}{r}\sigma^*\mO_{\mX(\vec{\lambda})}(1) )\big\}_{\vec{\lambda}\in (\bZ^+)^{\oplus l}}, $$
from the given weakly special collection as we vary $\vec{\lambda}\in (\bZ^+)^{\oplus l}$, where $\Delta_{\mX(\vec{\lambda})'}$ is the extension of $\Delta$ on $\mX(\vec{\lambda})'$. By Remark \ref{rem: crepant tc} we will see that all these 
$$\big\{(\mX(\vec{\lambda})', \Delta_{\mX(\vec{\lambda})'}\big\}_{\vec{\lambda}\in (\bZ^+)^{\oplus l}}$$ 
are crepant to each other. Let $\mX(\vec{\lambda})'_{0,i}$ be the component of $\mX(\vec{\lambda})'_0$ which corresponds to $c_i\cdot \ord_{E_i}$, and still write $c_i\cdot \ord_{E_{i,\infty}}$ for the valuation over $\mX(\vec{\lambda})'$ induced by $c_i\cdot \ord_{E_i}$. By the construction, we see 
$$\ord_{\mX(\vec{\lambda})'_{0,i}} =\lambda_i\cdot \ord_t+c_i\cdot \ord_{E_{i,\infty}},$$
where $t$ is the parameter of $\bA^1$. We also note here that, the data of the vector of positive integers somewhat gives the information on the polarization of the corresponding weakly special semi-test configuration. For example, if $\vec{\lambda}:=(\lambda_1,...,\lambda_l)$ and $\vec{\lambda'}:=(\lambda_1',...,\lambda_l')$ are two vectors satisfying that $\lambda_i-\lambda_i'=a$ is constant for $1\leq i\leq l$, then we could make sure that the two corresponding weakly special semi-test configurations share the same ambient space, and the two polarizations  are different by a multiple of the central fiber (the multiple is $a$).
\end{remark}

\begin{remark}\label{rem: crepant tc}
In this remark, we give an alternative way to construct weakly special degeneration\footnote{This means that we do not care about the information on the polarization.} from a weakly special collection, following the spirit of Proposition \ref{lc place to dege}. Given a weakly special collection
$\{c_1\cdot \ord_{E_1},...,c_l\cdot \ord_{E_l}\}$ and a vector $\vec{\lambda}:=(\lambda_1,...,\lambda_l)\in (\bZ^+)^{\oplus l}$, we first take a dlt modification $Y\to X$ which exactly extracts all $E_i$ (by the similar argument as \cite[Lemma 5.1]{CZ21} or Claim \ref{claim: extraction}).
Then we consider the product space 
$$Y\times \bA^1\to X\times \bA^1.$$
Let us write $X_0, Y_0$  for the central fibers and still denote $E_{i,\infty}$ to be the extension of $E_i$ on $Y\times \bA^1$. Thus we get a collection of divisorial valuations
$$\big\{ \lambda_j\cdot\ord_{Y_0}+ c_j\cdot \ord_{E_{j,\infty}}\mid j=1,...,l    \big\}$$
over $X\times \bA^1$. Let $\mE_j$ be a prime divisor over $X\times \bA^1$ such that
$$\ord_{\mE_j}:=\lambda_j\cdot\ord_{Y_0}+ c_j\cdot \ord_{E_{j,\infty}},$$
we see that all $\mE_j$ are log canonical places of $(X\times \bA^1, \Delta\times \bA^1+X_0)$ centered in $X_0$. Apply the similar argument as Claim \ref{claim: extraction}, one can find a $\bC^*$-equivariant birational model of $X\times \bA^1$ which exactly extracts all $\mE_j$, denoted by
$$\mY(\vec{\lambda})\to X\times \bA^1. $$
The central fiber of $\mY(\vec{\lambda})\to \bA^1$ clearly contains the trivial component, i.e., the strict transform of $X_0$. However, one can contract this trivial component via an MMP sequence as in the proof of Theorem \ref{thm: g-key}, denoted by 
$$\mY(\vec{\lambda})\dashrightarrow \mX(\vec{\lambda})'. $$
If we denote $\Delta_{\mY(\vec{\lambda})}$ (resp. $\Delta_{\mX(\vec{\lambda})'}$) to be the extension of $\Delta$ on $\mY(\vec{\lambda})$ (resp. $\mX(\vec{\lambda})'$), then it is not hard to see that the following three pairs are crepant to each other:
$$(\mY(\vec{\lambda}), \Delta_{\mY(\vec{\lambda})}), \ \ (\mX(\vec{\lambda})', \Delta_{\mX(\vec{\lambda})'}),  \ \  (X\times \bA^1, \Delta\times \bA^1),$$
and $(\mX(\vec{\lambda})', \Delta_{\mX(\vec{\lambda})'})$ here is exactly what we get in Remark \ref{rem: set of tc} up to isomorphism in codimension one.
\end{remark}

\

\section{Remark on global correspondence for log Fano pairs}

In the work \cite{CZ21}, for a given log canonical log Fano pair, we establish the correspondence between lc places of complements and weakly special test configurations with integral central fibers (see \cite[Theorem 1.2]{CZ21}). The method used there heavily depends on the computation of log canonical slope (see \cite[Theorem 3.1 and Theorem 6.2]{CZ21}), which follows the spirit of the proof of \cite[Theorem 4.10]{Xu21}. However, this method seems hard to apply to weakly special test configurations with reduced central fibers and to establish the correspondence as Theorem \ref{thm: main3} and \cite[Theorem A.2]{BLX19}. In this section, we give a remark on this topic.
Recall that a projective log pair $(X, B)$ is called \emph{log Fano} if $-K_X-B$ is ample. A log canonical log Fano pair means a log Fano pair with log canonical singularities. 

\begin{definition}
Let $(X, B)$ be an lc log Fano pair and $L:=-K_X-B$. A test configuration (resp. semi-test configuration) $(\mX, \mB;\mL)\to \bA^1$ of $(X, B;L)$ is \emph{weakly special} if $(\mX, \mB+\mX_0)$ is log canonical and $\mL\sim_\bQ-K_\mX-\mB$ is relatively ample (resp. semiample).  
\end{definition}

Let $(X, B)$ be an lc log Fano pair and $L:=-K_X-B$. Choose a positive integer $r$ such that $rL$ is Cartier. Let $(Z, B_Z)$ be the affine cone over $(X,B)$ with respect to $rL$, where $B_Z$ is the extension of $B$ on $Z$. Denote by $o\in Z$ the cone vertex, and we still write $v_0$ for the canonical valuation obtained by the blow-up of $o\in Z$.

\begin{definition}
Notation as above. A finite set of divisorial valuations over $X$, denoted by 
$$\{c_1\cdot \ord_{E_1},...,c_l\cdot\ord_{E_l}\},$$ 
is called a \emph{weakly special collection} of $(X,B)$ if the following two conditions are satisfied:
\begin{enumerate}
\item for each $i$,  $c_i\in \bZ^+$ and $E_i$ is a prime divisor over $X$;
\item there exists some effective $\bQ$-divisor $D\sim_\bQ -K_X-B$ such that $(X, B+D)$ is log canonical and all $E_i$ are lc places of $(X, B+D)$.
\end{enumerate}
\end{definition}

We also note that the definition of weakly special collection here is different from that in \cite[Appendix A]{BLX19}, as we already include the information of lc places. It is natural to ask the following question:

\begin{question}
Suppose that $(X,B)$ is an lc log Fano pair. Is there a correspondence 
between weakly special semi-test configurations of $(X, B;-K_X-B)$ and weakly special collections of $(X, B)$?
\end{question}

\bibliographystyle{amsalpha}
\bibliography{reference.bib}
\end{document}